\documentclass[preprint,12pt]{elsarticle}
\usepackage{amsthm,amsmath,amssymb}
\usepackage[colorlinks=true,citecolor=black,linkcolor=black,urlcolor=blue]{hyperref}
\usepackage{graphicx}
\usepackage{float}
\usepackage{mathtools}
\usepackage{epstopdf}
\usepackage{epsfig}
\usepackage{subfigure}
\usepackage{enumitem}
\usepackage{tikz}
\usepackage{nicematrix}
\usetikzlibrary{calc,decorations.pathmorphing,arrows.meta,fit,shapes.geometric}
\usepackage{arydshln}
\usepackage{blkarray}

\theoremstyle{plain}
\newtheorem{theorem}{Theorem}
\newtheorem{lemma}[theorem]{Lemma}
\newtheorem{corollary}[theorem]{Corollary}
\newtheorem{proposition}[theorem]{Proposition}

\theoremstyle{definition}
\newtheorem{definition}[theorem]{Definition}

\newtheorem{example}[theorem]{Example}

\newtheorem{problem}[theorem]{Problem}

\newtheoremstyle{myremark}
{3pt}   
{3pt}   
{\normalfont} 
{}      
{\bfseries} 
{.}     
{.5em}  
{}

\theoremstyle{myremark}
\newtheorem{remark}{Remark}

\newcommand{\GF}{\mathrm{GF}}

\title{Strict unimodality without log-concavity
for partial Petrial polynomials}

\author{Ruiqing Feng$^{1}$, Qi Yan$^{1}$\footnote{Corresponding author.}\\
\small School of Mathematics and Statistics, Lanzhou University, PR China\\
\small{\tt Email: fengrq2024@lzu.edu.cn; yanq@lzu.edu.cn}\\
}

\date{}
\journal{}

\newcommand{\rk}{\operatorname{rank}}
\newcommand{\crk}{\operatorname{corank}}
\newcommand{\PX}[1]{P_{#1}^{\times}}
\newcommand{\A}{\mathbf A}

\begin{document}
\begin{frontmatter}
\begin{abstract}
For a simple graph $G$, the partial Petrial
polynomial $P_G^{\times}(z)$ enumerates all diagonal
perturbations of its adjacency matrix over $\GF(2)$
by rank.
We derive a rooted attachment identity and apply it
to the tree $R_m$ obtained by joining a new vertex
to the middle vertex of each of $m$ disjoint copies
of the five-vertex path.
We prove that the nonzero coefficient sequence of
$P_{R_m}^{\times}(z)$ is strictly unimodal for every
$m\ge1$, but is log-concave if and only if
$1\le m\le3$.
Thus the $21$-vertex tree $R_4$ provides an explicit
strictly unimodal, non-log-concave example.
For $m\ge5$, the unique mode occurs at degree
$\lfloor(149m+85)/32\rfloor$.
More generally, for any fixed path rooted at an
internal vertex, joining a new vertex to the roots
of sufficiently many disjoint copies produces
a tree whose partial Petrial polynomial is not
log-concave.
\end{abstract}

\end{frontmatter}

\section{Introduction}

Unimodality and log-concavity are widely studied properties
of coefficient sequences of combinatorial polynomials.
Major results include Huh's proof of log-concavity for the
absolute values of chromatic polynomial
coefficients~\cite{Huh} and the corresponding theorem for
matroid characteristic polynomials established by
Adiprasito, Huh and Katz~\cite{AdiprasitoHuhKatz}.
For a nonnegative sequence with no internal zeros,
log-concavity implies unimodality.
Understanding when these properties hold, and when they
differ, is a natural problem for other graph polynomials.

The partial Petrial polynomial originated in the study
of ribbon graphs.
A partial Petrial is obtained by adding half-twists to
a specified subset of the edges of a ribbon graph.
Gross, Mansour and Tucker introduced the polynomial
that enumerates these operations according to
Euler genus~\cite{Gross}.
For bouquets, which are ribbon graphs with a single
vertex, this polynomial admits a description in terms
of binary matrices associated with intersection graphs.
Yan and Li showed that the partial Petrial polynomial
of a bouquet depends only on its intersection graph,
thereby defining the polynomial for circle
graphs~\cite{YanLi}.
Yu, Hao, Liu and Li obtained representations for
connected ribbon graphs using auxiliary
bouquets~\cite{YuHaoLiuLi}.
The matrix framework of Deng, Jin and
Yan~\cite{DengJinYan} extends the definition to
arbitrary simple graphs.

We work with this simple-graph formulation.
Let $G$ be a simple graph and let
$F\subseteq V(G)$.
Denote by $\A(G,F)$ the symmetric matrix over $\GF(2)$
whose off-diagonal entries agree with the adjacency
matrix of $G$ and whose diagonal entry at a vertex
is $1$ precisely when that vertex belongs to $F$.
The partial Petrial polynomial of $G$ is
\begin{equation}\label{eq:pgdef}
 \PX{G}(z)
 =\sum_{F\subseteq V(G)}z^{\rk\A(G,F)}.
\end{equation}
Thus its coefficients count diagonal perturbations
of the adjacency matrix by rank.
When $G$ is the intersection graph of a bouquet,
$\PX{G}(z)$ agrees with the topological partial Petrial
polynomial of that bouquet.
However, this correspondence does not identify
the polynomial in~\eqref{eq:pgdef} with the partial
Petrial polynomial of an arbitrary ribbon graph.
In particular, a coefficient example in the
ribbon-graph setting does not automatically provide
an example for~\eqref{eq:pgdef}.
Throughout this paper, partial Petrial polynomials
of graphs refer to the matrix-rank polynomials
defined in~\eqref{eq:pgdef}.

Several explicit formulas and structural results are
known for these polynomials.
Yan and Li determined the polynomials of complete
graphs and paths~\cite{YanLi}.
Feng, Yan and Zheng proved that a connected circle
graph has a binomial partial Petrial polynomial,
that is, a polynomial with exactly two nonzero terms,
if and only if it is a path~\cite{FengYanZheng}.
Deng and Zhang extended this characterization to
all connected simple graphs and characterized trees
with trinomial partial Petrial
polynomials~\cite{DengZhang}.
Feng, Guo and Yan established an edge recurrence
and proved that the lowest-degree nonzero coefficient
does not exceed the highest-degree coefficient,
with equality precisely for edgeless
graphs~\cite{FengGuoYan}.

These results lead to questions about the shape of
the entire coefficient sequence.
The polynomials in~\eqref{eq:pgdef} have no internal
zeros~\cite{DengJinYan}, so log-concavity implies
unimodality.
Paths give immediate examples of log-concavity,
since their polynomials have exactly two nonzero terms.
It is therefore natural to ask whether unimodality
can occur without log-concavity for these matrix-rank
polynomials, and whether such a distinction already
occurs for trees.
We answer both questions by constructing an explicit
infinite family whose coefficient sequences are
strictly unimodal but not log-concave.

Our starting point is a rooted attachment identity
valid for arbitrary simple graphs.
We apply it to the tree $R_m$ obtained by joining
a new vertex to the middle vertex of each of $m$
vertex-disjoint copies of the five-vertex path.
Thus $|V(R_m)|=5m+1$.
We prove that the nonzero coefficient sequence of
$\PX{R_m}(z)$ is strictly unimodal for every $m\ge1$,
whereas it is log-concave if and only if
$1\le m\le3$
(Theorems~\ref{st:3.5} and~\ref{st:3.8}).
For $m\ge5$, its unique mode occurs at degree
\[
 \mu_m=\left\lfloor\frac{149m+85}{32}\right\rfloor.
\]
Consequently, $\{R_m:m\ge4\}$ is an infinite family
showing that strict unimodality does not imply
log-concavity for the partial Petrial polynomials
of simple graphs, even when the graphs are trees.
The member $R_4$ provides an explicit example on
$21$ vertices.

The attachment identity also yields a sufficient
criterion for eventual failure of log-concavity
under repeated attachment.
Applying this criterion to paths, we prove that
for every fixed path rooted at an internal vertex,
joining a new vertex to the roots of sufficiently
many vertex-disjoint copies produces a tree whose
partial Petrial polynomial is not log-concave
(Corollary~\ref{st:3.12}).
Thus the failure of log-concavity in the family
$R_m$ is part of a more general phenomenon
associated with repeated path branches.

The paper is organized as follows.
Section~2 introduces the corank polynomial,
coefficient conventions and leaf reduction.
Section~3 establishes the rooted attachment identity,
constructs $R_m$ and determines its exact
log-concavity threshold.
Section~4 proves strict unimodality and determines
the mode locations.
Section~5 develops the criterion for eventual failure
of log-concavity and applies it to repeated path
branches.
Section~6 concludes with questions on unimodality
and structural conditions for log-concavity.

\section{Preliminaries}

For vertex-disjoint graphs $G_1$ and $G_2$, we write
$G_1\sqcup G_2$ for their disjoint union.
For a simple graph $G$ with $n=|V(G)|$, define
\[
 \crk\A(G,F)=n-\rk\A(G,F),
 \qquad
 C_G(t)=\sum_{F\subseteq V(G)}t^{\crk\A(G,F)}.
\]
By~\eqref{eq:pgdef},
\begin{equation}\label{eq:reverse}
 \PX{G}(z)=z^nC_G(z^{-1}).
\end{equation}
Thus the nonzero coefficient sequences of $\PX{G}(z)$
and $C_G(t)$ are reversals of one another.

Let
\[
 f(z)=\sum_{k=r}^{d}a_kz^k,
 \qquad a_r,a_d>0,\quad a_k\ge0.
\]
Unless otherwise indicated, its coefficient sequence
means $(a_r,\ldots,a_d)$, including all coefficients
between the lowest and highest nonzero degrees.
The polynomial is called \emph{interpolating} if
$a_k>0$ for every $r\le k\le d$.
Partial Petrial polynomials have this
property~\cite{DengJinYan}.

The coefficient sequence is \emph{log-concave} if
\[
 a_k^2\ge a_{k-1}a_{k+1}
 \qquad(r<k<d),
\]
and \emph{strictly log-concave} if all these inequalities
are strict.
It is \emph{unimodal} if it is first nondecreasing and
then nonincreasing, and \emph{strictly unimodal} if it
increases strictly up to a unique mode and decreases
strictly thereafter.
The mode is allowed to occur at either endpoint.
For nonnegative sequences with no internal zeros,
log-concavity implies unimodality.

All four properties are invariant under reversal.
Consequently, we may study either $\PX{G}(z)$ or $C_G(t)$,
as convenient.

We will use the following leaf-reduction identity
to derive the path formula needed later.

\begin{lemma}[\cite{DengJinYan}]
\label{st:2.4}
Let $u$ be a leaf of a simple graph $G$ with
unique neighbor $v$.
Write $G-u$ and $G-\{u,v\}$ for the induced subgraphs
obtained by deleting the indicated vertices.
Then
\begin{equation}\label{eq:leaf}
 \PX{G}(z)
 =z\PX{G-u}(z)+2z^2\PX{G-\{u,v\}}(z).
\end{equation}
\end{lemma}

\section{Repeated-branch trees and the log-concavity threshold}

We construct a family of trees $R_m$ and determine
exactly when their partial Petrial polynomials are
log-concave.
We begin with an attachment formula for arbitrary
rooted simple graphs, derived directly from kernel
equations over $\GF(2)$.

Let $H$ be a simple graph with a distinguished vertex
$c\in V(H)$. For $m\ge1$, take $m$ vertex-disjoint copies
$H_1,\ldots,H_m$ of $H$, with corresponding distinguished
vertices $c_1,\ldots,c_m$.
Add a new vertex $r$ and the edges $rc_i$ for $1\le i\le m$,
and denote the resulting graph by $T_m(H,c)$.
If $H$ is a tree, then so is $T_m(H,c)$.
The family $R_m$ studied below corresponds to $H=P_5$,
with $c$ the middle vertex.

To determine the corank polynomial of $T_m(H,c)$, we need more than the corank of each branch perturbation: we must also know whether every vector in its kernel vanishes at the attachment vertex. Accordingly, define
\[
 U_{H,c}(t)=
 \sum_{\substack{F\subseteq V(H)\\
 x_c=0\ \text{for every }\mathbf x\in\ker\A(H,F)}}
 t^{\crk\A(H,F)}.
\]
All nonsingular perturbations are included in $U_{H,c}(t)$, since their kernels are trivial.
The complementary polynomial $C_H(t)-U_{H,c}(t)$ counts the perturbations whose kernels contain a vector with coordinate $1$ at $c$.

The attachment formula below distinguishes two cases:
either every branch perturbation is counted by
$U_{H,c}(t)$, or at least one is counted by
$C_H(t)-U_{H,c}(t)$.

\begin{lemma}\label{st:root-general}
For every rooted simple graph $(H,c)$ and every $m\ge1$,
\begin{equation}\label{eq:root-general}
 C_{T_m(H,c)}(t)
 =(1+t)U_{H,c}(t)^m
 +\frac2t\bigl[C_H(t)^m-U_{H,c}(t)^m\bigr].
\end{equation}
\end{lemma}

\begin{proof}
All matrix and vector calculations below are over $\GF(2)$,
while the generating functions have integer coefficients.
Write $h=|V(H)|$.
Fix a diagonal perturbation $F_i\subseteq V(H_i)$ on each
copy $H_i$, and let
\[
 B_i=\A(H_i,F_i),\qquad
 \nu_i=\dim\ker B_i,\qquad
 K=\sum_{i=1}^m\nu_i.
\]
The diagonal entry $\delta_r\in\GF(2)$ at the new root
remains to be chosen.

Let $\mathbf e_i\in\GF(2)^h$ be the standard basis vector
corresponding to $c_i$.
A vector with coordinate $s\in\GF(2)$ at $r$ and branch
components $\mathbf u_i\in\GF(2)^h$ lies in the kernel
of the full perturbed matrix if and only if it satisfies
the branch equations
\[
 B_i\mathbf u_i+s\mathbf e_i=\mathbf0
 \qquad(1\le i\le m)
\]
and the root equation
\[
 \delta_rs+\sum_{i=1}^m
 \mathbf e_i^{\mathsf T}\mathbf u_i=0.
\]
We determine the corank for each of the two choices
of $\delta_r$.

\medskip
\noindent
\emph{Case 1: Every chosen branch perturbation is counted by
$U_{H,c}(t)$.}

In this case,
\[
 \mathbf e_i^{\mathsf T}\mathbf z_i=0
 \qquad\text{for every }\mathbf z_i\in\ker B_i
 \quad(1\le i\le m).
\]
Since $B_i$ is symmetric,
\[
 \operatorname{im}B_i=(\ker B_i)^\perp.
\]
Consequently, we may choose a vector $\mathbf v_i$
satisfying $B_i\mathbf v_i=\mathbf e_i$.

For each fixed $s\in\GF(2)$, the branch equations have
the general solution
\[
 \mathbf u_i=s\mathbf v_i+\mathbf z_i,
 \qquad \mathbf z_i\in\ker B_i
 \quad(1\le i\le m).
\]
The vectors $\mathbf z_i$ vary independently, giving
$K$ free parameters.
Moreover,
\[
 \mathbf e_i^{\mathsf T}\mathbf u_i
 =s\,\mathbf e_i^{\mathsf T}\mathbf v_i,
\]
because $\mathbf e_i^{\mathsf T}\mathbf z_i=0$.
Thus none of these $K$ parameters appears in the root equation.
Putting
\[
 \lambda=\sum_{i=1}^m
 \mathbf e_i^{\mathsf T}\mathbf v_i\in\GF(2),
\]
the root equation becomes
\[
 (\delta_r+\lambda)s=0.
\]

If $\delta_r=\lambda$, then $s$ is an additional free parameter,
and the corank is $K+1$.
For the other choice of $\delta_r$, the equation forces $s=0$,
and the corank is $K$.
Hence, for this fixed tuple of branch perturbations,
the two diagonal choices at the root contribute
\[
 t^K+t^{K+1}=(1+t)t^K.
\]

\medskip
\noindent
\emph{Case 2: At least one chosen branch perturbation is not
counted by $U_{H,c}(t)$.}

There is an index $j$ and a vector
$\boldsymbol\eta_j\in\ker B_j$ such that
\[
 \mathbf e_j^{\mathsf T}\boldsymbol\eta_j=1.
\]
Since $B_j$ is symmetric and
$B_j\boldsymbol\eta_j=\mathbf0$, multiplying the $j$th
branch equation on the left by
$\boldsymbol\eta_j^{\mathsf T}$ gives
\[
 0
 =\boldsymbol\eta_j^{\mathsf T}
   (B_j\mathbf u_j+s\mathbf e_j)
 =s.
\]
Thus $s=0$, and every branch equation reduces to
$\mathbf u_i\in\ker B_i$.

The branch vectors therefore range over the $K$-dimensional
space
\[
 W=\ker B_1\times\cdots\times\ker B_m.
\]
On this space, the root equation is
\[
 L(\mathbf u_1,\ldots,\mathbf u_m)
 :=\sum_{i=1}^m\mathbf e_i^{\mathsf T}\mathbf u_i=0.
\]
The linear functional $L$ is nonzero: it takes the value $1$
on the tuple with $\boldsymbol\eta_j$ in the $j$th component
and zero in every other component.
Since $L:W\to\GF(2)$ is nonzero, it has rank $1$.
The rank--nullity theorem therefore gives
\[
 \dim\ker L=\dim W-1=K-1.
\]

Since $s=0$, the value of $\delta_r$ does not affect any
of these equations.
Both diagonal choices at the root therefore give corank
$K-1$, and their combined contribution is
\(
 2t^{K-1}.
\)

\medskip
It remains to sum over all tuples $(F_1,\ldots,F_m)$.
Each tuple has weight
\[
 \prod_{i=1}^m t^{\nu_i}=t^K
\]
in $C_H(t)^m$.
The tuples in Case~1 are counted by $U_{H,c}(t)^m$,
and each weight $t^K$ is replaced by $(1+t)t^K$.
The remaining tuples are counted by
$C_H(t)^m-U_{H,c}(t)^m$,
and each weight $t^K$ is replaced by $2t^{K-1}$.
Summing these contributions proves~\eqref{eq:root-general}.

Finally, every tuple in Case~2 contains a branch with a
nonzero kernel vector, so $K\ge1$.
Hence $C_H(t)^m-U_{H,c}(t)^m$ is divisible by $t$,
and the right-hand side of~\eqref{eq:root-general}
is a polynomial.
\end{proof}

We now specialize the general attachment identity to the case
in which
\[
 \crk\A(H,F)\le1
 \qquad\text{for every }F\subseteq V(H).
\]
Under this assumption, each singular matrix $\A(H,F)$
has a one-dimensional kernel and hence a unique nonzero
kernel vector over $\GF(2)$.
Denote this vector by $\boldsymbol\kappa(F)$ and its
coordinate at $c$ by $\kappa_c(F)$.

We classify the diagonal perturbations into three types,
with counts
\[
 \begin{aligned}
 a&=\#\{F\subseteq V(H):
       \crk\A(H,F)=0\},\\
 b&=\#\{F\subseteq V(H):
       \crk\A(H,F)=1,\ \kappa_c(F)=0\},\\
 \gamma&=\#\{F\subseteq V(H):
       \crk\A(H,F)=1,\ \kappa_c(F)=1\}.
 \end{aligned}
\]
These three classes partition all diagonal perturbations, so
\[
 a+b+\gamma=2^{|V(H)|}.
\]

The polynomial $U_{H,c}(t)$ counts all nonsingular
perturbations and precisely those singular perturbations
for which $\kappa_c(F)=0$.
Consequently,
\[
 C_H(t)=a+(b+\gamma)t,
 \qquad
 U_{H,c}(t)=a+bt.
\]

\begin{corollary}\label{st:3.1}
Let $(H,c)$ be a rooted simple graph satisfying
$\crk\A(H,F)\le1$ for every $F\subseteq V(H)$.
With $a$, $b$ and $\gamma$ as defined above, we have,
for every $m\ge1$,
\begin{equation}\label{eq:rootjoin}
 C_{T_m(H,c)}(t)
 =(1+t)(a+bt)^m
 +\frac2t\bigl[
   (a+(b+\gamma)t)^m-(a+bt)^m
 \bigr].
\end{equation}
\end{corollary}

\begin{proof}
Substituting $U_{H,c}(t)=a+bt$ and
$C_H(t)=a+(b+\gamma)t$ into~\eqref{eq:root-general}
gives~\eqref{eq:rootjoin}.
\end{proof}

\begin{figure}[htbp]
 \centering
 \begin{tikzpicture}[
   x=1cm,
   y=1cm,
   scale=.84,
   every node/.style={font=\scriptsize}
 ]
  \node[
    circle,draw,inner sep=1.2pt,label=above:$r$
  ] (root) at (0,3.1) {};

  \foreach \i/\xc in {1/-5.4,2/-1.8,3/1.8,4/5.4}{
   \node[
     circle,draw,inner sep=1.2pt,label=above:$c_{\i}$
   ] (c\i) at (\xc,1.7) {};

   \node[
     circle,draw,inner sep=1.2pt,
     label=below left:$u_{\i,1}$
   ] (x\i) at (\xc-.7,.6) {};

   \node[
     circle,draw,inner sep=1.2pt,
     label=below right:$u_{\i,2}$
   ] (y\i) at (\xc+.7,.6) {};

   \node[
     circle,draw,inner sep=1.2pt,
     label=below left:$v_{\i,1}$
   ] (xx\i) at (\xc-.7,-.55) {};

   \node[
     circle,draw,inner sep=1.2pt,
     label=below right:$v_{\i,2}$
   ] (yy\i) at (\xc+.7,-.55) {};

   \draw
     (root)--(c\i)
     (c\i)--(x\i)
     (c\i)--(y\i)
     (x\i)--(xx\i)
     (y\i)--(yy\i);
  }
 \end{tikzpicture}
 \caption{The tree $R_4$, formed by joining $r$ to the middle
 vertex of each of four copies of $P_5$.}
 \label{fig:R4}
\end{figure}
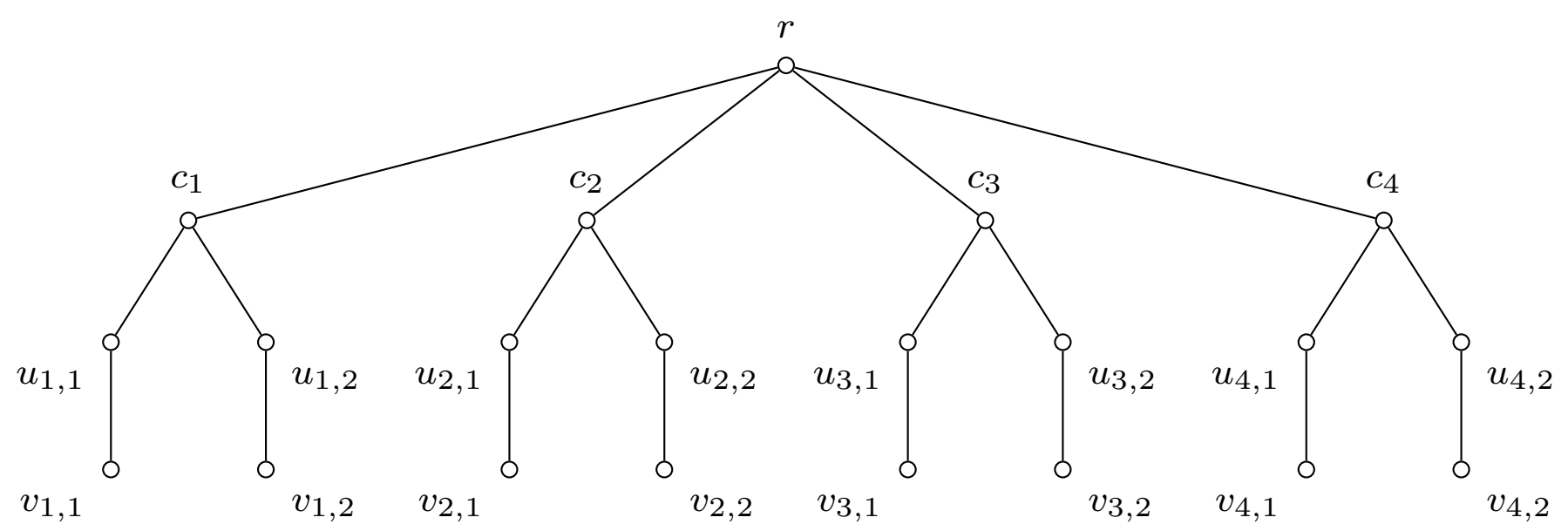

\begin{definition}\label{st:3.2}
Let $P_5$ be a five-vertex path rooted at its middle
vertex $c$, and set $R_m=T_m(P_5,c)$ for $m\ge1$.
Thus $R_m$ is a tree on $5m+1$ vertices.
Figure~\ref{fig:R4} depicts $R_4$.
\end{definition}

\begin{proposition}\label{st:3.3}
For $P_5$ rooted at its middle vertex $c$, every diagonal
perturbation has corank at most $1$, and the counts
defined above satisfy
\[
 a=21,\qquad b=2,\qquad \gamma=9.
\]
Consequently,
\[
 C_{P_5}(t)=21+11t,
 \qquad
 U_{P_5,c}(t)=21+2t.
\]
\end{proposition}

\begin{proof}
All matrix and vector calculations below are over $\GF(2)$.
Label the vertices of the path, in path order, as
$v_1,u_1,c,u_2,v_2$.
Using the matrix order $c,u_1,v_1,u_2,v_2$, write
\[
 A=\A(P_5,F)=
 \begin{pmatrix}
 \delta_c&1&0&1&0\\
 1&\alpha_1&1&0&0\\
 0&1&\beta_1&0&0\\
 1&0&0&\alpha_2&1\\
 0&0&0&1&\beta_2
 \end{pmatrix},
\]
where the five diagonal entries range independently
over $\GF(2)$.

Let
\[
 \mathbf x=(s,p_1,q_1,p_2,q_2)^{\mathsf T}\in\ker A.
\]
The equations at the two vertices of arm $i$ are
\[
 \alpha_i p_i+q_i+s=0,\qquad
 p_i+\beta_iq_i=0
 \quad(i=1,2),
\]
and the equation at the central vertex is
\[
 \delta_cs+p_1+p_2=0.
\]
For each $i$, the second arm equation gives $p_i=\beta_iq_i$.
Substituting this into the first gives
\[
 (1+\alpha_i\beta_i)q_i=s.
\]

The block corresponding to the vertices $u_i,v_i$ is
\[
 D_i=
 \begin{pmatrix}
 \alpha_i&1\\
 1&\beta_i
 \end{pmatrix},
\]
with determinant $1+\alpha_i\beta_i$.
Thus it is nonsingular for the three choices
\[
 (\alpha_i,\beta_i)=(0,0),(0,1),(1,0),
\]
and singular only for $(\alpha_i,\beta_i)=(1,1)$.
We distinguish the following three cases.

\medskip
\noindent
\emph{Case 1: Both arm blocks are nonsingular.}

Here $1+\alpha_i\beta_i=1$ for $i=1,2$, so
\[
 q_i=s,\qquad p_i=\beta_i s.
\]
Consequently, every kernel vector has the form
\[
 \mathbf x
 =s(1,\beta_1,1,\beta_2,1)^{\mathsf T},
\]
and the central equation becomes
\[
 (\delta_c+\beta_1+\beta_2)s=0.
\]

For each fixed choice of the two arm blocks, exactly one
value of $\delta_c$, namely $\delta_c=\beta_1+\beta_2$,
leaves $s$ free.
The resulting kernel is one-dimensional, and its nonzero
vector has central coordinate $1$.
The other value of $\delta_c$ forces $s=0$, so the full
matrix is nonsingular.

There are $3^2=9$ choices of the two nonsingular arm blocks.
This case therefore contributes $9$ to $\gamma$ and $9$ to $a$.

\medskip
\noindent
\emph{Case 2: Exactly one arm block is singular.}

By symmetry, assume that the first arm block is singular.
Then $\alpha_1=\beta_1=1$, so the arm equations force
$s=0$ and $p_1=q_1$.
Since the second arm block is nonsingular, we obtain
$q_2=p_2=0$.
The central equation now gives $p_1=0$.
Hence every coordinate of $\mathbf x$ is zero, and the full
matrix is nonsingular for both choices of $\delta_c$.

There are two choices for which arm is singular,
three choices for the nonsingular arm block, and two choices
for $\delta_c$.
This case therefore contributes $2\cdot3\cdot2=12$ to $a$.

\medskip
\noindent
\emph{Case 3: Both arm blocks are singular.}

Now $\alpha_i=\beta_i=1$ for $i=1,2$.
The arm equations give
\[
 s=0,\qquad p_1=q_1,\qquad p_2=q_2.
\]
The central equation reduces to $p_1+p_2=0$.
Therefore
\[
 \ker A
 =
 \left\{
 \tau(0,1,1,1,1)^{\mathsf T}:
 \tau\in\GF(2)
 \right\}.
\]
This kernel is one-dimensional, and its nonzero vector
has central coordinate $0$.
Both choices of $\delta_c$ give this kernel, so this case
contributes $2$ to $b$.

\medskip
These three cases exhaust all diagonal perturbations and
show that every perturbed matrix has corank at most $1$.
Adding their contributions gives
\[
 a=9+12=21,\qquad b=2,\qquad \gamma=9.
\]
The stated formulas for $C_{P_5}(t)$ and $U_{P_5,c}(t)$
follow from the preceding classification.
\end{proof}

\begin{proposition}\label{st:3.4}
Put $C_m(t)=C_{R_m}(t)$ for $m\ge1$.
Then
\begin{equation}\label{eq:cm}
 C_m(t)
 =(1+t)(21+2t)^m
 +\frac2t\bigl[(21+11t)^m-(21+2t)^m\bigr].
\end{equation}
Moreover, setting $C_0(t)=1+t$, we have, for every $m\ge0$,
\begin{equation}\label{eq:cmrecur}
 C_{m+1}(t)
 =(21+2t)C_m(t)+18(21+11t)^m.
\end{equation}
\end{proposition}

\begin{proof}
By Proposition~\ref{st:3.3}, every diagonal perturbation
of $P_5$ has corank at most $1$, and
\[
 U_{P_5,c}(t)=21+2t,
 \qquad
 C_{P_5}(t)=21+11t.
\]
Thus Corollary~\ref{st:3.1} gives~\eqref{eq:cm}.

For the recurrence, put
\[
 p=21+2t,\qquad q=21+11t.
\]
Formula~\eqref{eq:cm} also holds at $m=0$, since its
right-hand side then equals $1+t=C_0(t)$.
Hence, for every $m\ge0$,
\[
 \begin{aligned}
 C_{m+1}(t)-pC_m(t)
 &=\frac2t\bigl(q^{m+1}-pq^m\bigr)\\
 &=\frac2t q^m(q-p)\\
 &=18q^m,
 \end{aligned}
\]
where the last equality uses $q-p=9t$.
Substituting back for $p$ and $q$
gives~\eqref{eq:cmrecur}.
\end{proof}

For $m\ge1$, the difference-of-powers identity
\[
 x^m-y^m=(x-y)\sum_{j=0}^{m-1}x^{m-1-j}y^j
\]
allows us to rewrite~\eqref{eq:cm} without subtraction.
Applying this identity with $x=21+11t$ and $y=21+2t$,
and using $x-y=9t$, gives
\begin{equation}\label{eq:cm-sum}
 C_m(t)
 =(1+t)(21+2t)^m
 +18\sum_{j=0}^{m-1}
 (21+11t)^{m-1-j}(21+2t)^j.
\end{equation}
The first term has strictly positive coefficients in every
degree from $0$ to $m+1$, while the remaining sum has
nonnegative coefficients and degree at most $m-1$.
Thus $C_m(t)$ has degree $m+1$ and strictly positive
coefficients  in every degree from $0$ to $m+1$.

To recover the rank polynomial, recall that $R_m$ has
$5m+1$ vertices.
Applying~\eqref{eq:reverse} to~\eqref{eq:cm} gives
\begin{equation}\label{eq:prm}
 \begin{aligned}
 \PX{R_m}(z)
 &=z^{5m+1}C_m(z^{-1})\\
 &=z^{4m}(1+z)(21z+2)^m\\
 &\quad
 +2z^{4m+2}
 \bigl[(21z+11)^m-(21z+2)^m\bigr].
 \end{aligned}
\end{equation}
In particular, the nonzero coefficients of $\PX{R_m}(z)$
occur precisely in degrees $4m,\ldots,5m+1$.

Finally, evaluating~\eqref{eq:cm} at $t=1$ gives
\[
 C_m(1)
 =2\cdot23^m+2(32^m-23^m)
 =2^{5m+1}.
\]
This agrees with the total number of diagonal perturbations
on the $5m+1$ vertices of $R_m$.

\begin{example}\label{ex:R4}
The tree $R_4$ has $5\cdot4+1=21$ vertices.
Substituting $m=4$ into~\eqref{eq:prm} gives
\begin{equation}\label{eq:R4-example}
 \begin{aligned}
 \PX{R_4}(z)={}&16z^{16}+688z^{17}+40506z^{18}\\
 &+306936z^{19}+887733z^{20}+861273z^{21}.
 \end{aligned}
\end{equation}
Its nonzero coefficients satisfy
\[
 16<688<40506<306936<887733>861273.
\]
Hence the nonzero coefficient sequence is strictly unimodal,
with its unique mode at degree $20$.

However,
\[
 688^2=473344<648096=16\cdot40506,
\]
so log-concavity fails at degree $17$.
Thus the nonzero coefficient sequence is not log-concave.
\end{example}

We now determine the exact log-concavity threshold
for the family $R_m$.

\begin{theorem}\label{st:3.5}
For $m\ge1$, the nonzero coefficient sequence of
$\PX{R_m}(z)$ is log-concave if and only if $1\le m\le3$.
\end{theorem}

\begin{proof}
We first consider $m=1,2,3$.
By~\eqref{eq:prm}, the corresponding nonzero coefficient
sequences are
\[
 \begin{array}{c|l}
 m & \text{Nonzero coefficient sequence}\\ \hline
 1 & (2,23,39)\\
 2 & (4,88,759,1197)\\
 3 & (8,260,5544,26649,33075).
 \end{array}
\]
For each sequence $(a_0,\ldots,a_r)$ above, the differences
$a_i^2-a_{i-1}a_{i+1}$ at the internal indices are
\[
 \begin{array}{c|l}
 m & \text{Internal log-concavity differences}\\ \hline
 1 & 451\\
 2 & 4708,\ 470745\\
 3 & 23248,\ 23807196,\ 526801401.
 \end{array}
\]
All these differences are positive.
Thus the nonzero coefficient sequences are strictly
log-concave for $1\le m\le3$.

Now let $m\ge4$ and write
\[
 C_m(t)=\sum_{k=0}^{m+1}d_{m,k}t^k.
\]
The nonzero coefficient sequences of $C_m(t)$ and
$\PX{R_m}(z)$ are reversals of one another, so it suffices
to show that
\[
 d_{m,m}^2<d_{m,m-1}d_{m,m+1}.
\]
This corresponds to failure of log-concavity at degree
$4m+1$ in $\PX{R_m}(z)$.

The second term on the right-hand side of~\eqref{eq:cm}
has degree $m-1$ and leading coefficient $2(11^m-2^m)$.
Applying the binomial theorem to the first term therefore gives
\[
 \begin{aligned}
 d_{m,m+1}
 &=2^m,\\
 d_{m,m}
 &=2^m+21m2^{m-1}
 =2^{m-1}(21m+2),\\
 d_{m,m-1}
 &=21m2^{m-1}
 +441\binom m2 2^{m-2}
 +2(11^m-2^m).
 \end{aligned}
\]
Set
\[
 \Delta_m=d_{m,m}^2-d_{m,m-1}d_{m,m+1},
 \qquad
 \Phi_m=441m^2+525m+24.
\]
Substituting these coefficients and simplifying gives
\begin{equation}\label{eq:lm}
 \Delta_m
 =2^{2m-3}
 \left[
   \Phi_m-16\left(\frac{11}{2}\right)^m
 \right].
\end{equation}
Since $2^{2m-3}>0$, it remains to prove
\begin{equation}\label{eq:4threshold}
 16\left(\frac{11}{2}\right)^m>\Phi_m
 \qquad(m\ge4).
\end{equation}

We prove this inequality by induction.
For $m=4$,
\[
 16\left(\frac{11}{2}\right)^4
 =14641>9180=\Phi_4.
\]
For the induction step, observe that
\[
 \frac{\Phi_m}{m^2}
 =441+\frac{525}{m}+\frac{24}{m^2}
\]
decreases as $m$ increases.
Hence, for $m\ge4$,
\[
 \Phi_{m+1}
 \le\left(1+\frac1m\right)^2\Phi_m
 \le\frac{25}{16}\Phi_m
 <\frac{11}{2}\Phi_m.
\]
Thus, if~\eqref{eq:4threshold} holds at $m$, then
\[
 16\left(\frac{11}{2}\right)^{m+1}
 >\frac{11}{2}\Phi_m
 >\Phi_{m+1}.
\]
Together with the base case, this proves
\eqref{eq:4threshold} for every $m\ge4$.

It follows from~\eqref{eq:lm} that $\Delta_m<0$
for every $m\ge4$.
Consequently, $C_m(t)$ violates the log-concavity inequality
at $d_{m,m}$.
By reversal, the nonzero coefficient sequence of
$\PX{R_m}(z)$ is also not log-concave.
Together with the cases $m=1,2,3$, this proves the theorem.
\end{proof}
	
Thus, for every $m\ge4$, log-concavity already fails at
degree $4m+1$, the second-lowest nonzero degree of
$\PX{R_m}(z)$.
	
\section{Strict unimodality and mode locations for $R_m$}

Theorem~\ref{st:3.5} shows that the nonzero coefficient
sequence of $\PX{R_m}(z)$ is not log-concave for $m\ge4$.
We now prove that it is strictly unimodal for every $m\ge1$
and determine its exact mode for $m\ge5$.
The cases $1\le m\le4$ are already determined by the
explicit coefficient sequences above.

Continue to write
\[
 C_m(t)=\sum_{k=0}^{m+1}d_{m,k}t^k.
\]
Since
\[
 \PX{R_m}(z)=z^{5m+1}C_m(z^{-1}),
\]
the two nonzero coefficient sequences are reversals of
one another.
Thus it suffices to study strict unimodality of $C_m(t)$.
Moreover, a unique mode of $C_m(t)$ at degree $q$
corresponds to a unique mode of $\PX{R_m}(z)$ at degree
$5m+1-q$.

We first record an explicit formula for $d_{m,k}$ and a
recurrence relating the coefficients for consecutive
values of $m$.

\begin{lemma}
\label{st:3.6}
For $m\ge1$ and $0\le k\le m$, define
\[
 u_{m,k}=\binom mk21^{m-k}2^k,
 \qquad
 b_{m,k}=\binom mk21^{m-k}11^k.
\]
Set $u_{m,k}=b_{m,k}=0$ when $k<0$ or $k>m$.
Also set $d_{m,k}=0$ when $k<0$ or $k>m+1$.
Then, for $0\le k\le m+1$,
\begin{equation}\label{eq:explicit}
 d_{m,k}
 =u_{m,k}+u_{m,k-1}
 +2\bigl(b_{m,k+1}-u_{m,k+1}\bigr).
\end{equation}
Moreover, for $0\le k\le m+2$,
\begin{equation}\label{eq:dmrecur}
 d_{m+1,k}
 =21d_{m,k}+2d_{m,k-1}+18b_{m,k}.
\end{equation}
\end{lemma}

\begin{proof}
By the binomial theorem, $u_{m,k}$ and $b_{m,k}$ are
the coefficients of $t^k$ in $(21+2t)^m$ and
$(21+11t)^m$, respectively.

Since $b_{m,0}=u_{m,0}=21^m$, the second term on the
right-hand side of~\eqref{eq:cm} is a polynomial whose
coefficient of $t^k$ is
\[
 2\bigl(b_{m,k+1}-u_{m,k+1}\bigr).
\]
The first term contributes $u_{m,k}+u_{m,k-1}$,
which proves~\eqref{eq:explicit}.

Similarly, comparing coefficients of $t^k$
in~\eqref{eq:cmrecur} gives
\[
 d_{m+1,k}
 =21d_{m,k}+2d_{m,k-1}+18b_{m,k},
\]
proving~\eqref{eq:dmrecur}.
The zero conventions ensure that both formulas hold
at the boundary indices as well.
\end{proof}

We first motivate a candidate for the mode of $C_m(t)$.
By~\eqref{eq:explicit}, the coefficient $d_{m,k}$ contains
the auxiliary term $2b_{m,k+1}$.
We therefore examine the shifted auxiliary sequence
$(b_{m,k+1})_k$.

From the definition of $b_{m,k}$,
\[
 \frac{b_{m,k+1}}{b_{m,k}}
 =\frac{\binom m{k+1}}{\binom mk}\frac{11}{21}
 =\frac{11(m-k)}{21(k+1)}
 \qquad(0\le k\le m),
\]
where the endpoint $k=m$ uses the zero convention.
Consequently,
\[
 \begin{aligned}
 b_{m,k+2}-b_{m,k+1}
 &=b_{m,k+1}
   \left(\frac{11(m-k-1)}{21(k+2)}-1\right)\\
 &=\frac{b_{m,k+1}}{21(k+2)}
   (11m-32k-53)
 \qquad(0\le k\le m-1).
 \end{aligned}
\]
The factor outside the final parentheses is positive.
Within this range, the difference is therefore positive
when $k<(11m-53)/32$ and negative when
$k>(11m-53)/32$.

For $m\ge5$, this suggests the candidate mode
\begin{equation}\label{eq:qm}
 q_m=\left\lceil\frac{11m-53}{32}\right\rceil.
\end{equation}
Indeed, $q_m$ is the first integer index at which the
auxiliary difference above is nonpositive.
If $(11m-53)/32$ is an integer, the shifted auxiliary
sequence has equal adjacent maxima at indices $q_m$
and $q_m+1$.
To establish $q_m$ as the unique mode of $C_m(t)$,
we must control the remaining terms in~\eqref{eq:explicit},
including their effect when the auxiliary sequence has
two equal maxima.

The candidate mode moves by at most one position as $m$
increases.
Since
\[
 \frac{11(m+1)-53}{32}
 =\frac{11m-53}{32}+\frac{11}{32},
 \qquad 0<\frac{11}{32}<1,
\]
we have
\[
 q_{m+1}-q_m\in\{0,1\}.
\]

To see how the auxiliary differences enter the induction,
take consecutive differences in~\eqref{eq:dmrecur}.
For $0\le k\le m+1$, this gives
\[
 \begin{aligned}
 d_{m+1,k+1}-d_{m+1,k}
 ={}&21(d_{m,k+1}-d_{m,k})\\
 &+2(d_{m,k}-d_{m,k-1})\\
 &+18(b_{m,k+1}-b_{m,k}).
 \end{aligned}
\]
In the induction below, the induction hypothesis,
together with this recurrence and the signs of the
auxiliary differences, will show that the coefficient
sequence of $C_{m+1}(t)$ increases strictly up to index
$q_m$ and decreases strictly from index $q_m+1$ onward.
The remaining comparison is therefore between
$d_{m+1,q_m}$ and $d_{m+1,q_m+1}$.
Since $q_{m+1}\in\{q_m,q_m+1\}$, this comparison will
determine whether the unique mode is $q_m$ or $q_m+1$.
We first establish the estimate needed for this comparison,
and then complete the induction.

\begin{lemma}
\label{st:3.7}
Let $m\ge5$ and $q=q_m$. Then
\[
 \begin{aligned}
 d_{m+1,q+1}&>d_{m+1,q}
 &&\text{if }q_{m+1}=q+1,\\
 d_{m+1,q+1}&<d_{m+1,q}
 &&\text{if }q_{m+1}=q.
 \end{aligned}
\]
\end{lemma}

\begin{proof}
Let
\[
 \Gamma_m=d_{m+1,q+1}-d_{m+1,q}.
\]
The definition of $q_m$ gives $1\le q\le m-2$.
Put
\[
 N_m=11(m+1)-32q-53,
 \qquad
 E_q=\left(\frac{11}{2}\right)^{q+1}-1.
\]
Since
\[
 \frac{11(m+1)-53}{32}=q+\frac{N_m}{32}
\]
and $q_{m+1}\in\{q,q+1\}$, we have
\begin{equation}\label{eq:qmstep}
 \begin{aligned}
 N_m>0&\iff q_{m+1}=q+1,\\
 N_m\le0&\iff q_{m+1}=q.
 \end{aligned}
\end{equation}
It therefore suffices to prove that $\Gamma_m>0$ when
$N_m>0$, and $\Gamma_m<0$ when $N_m\le0$.

We first derive an expression for this difference.
Subtracting~\eqref{eq:explicit}, with $m$ replaced by $m+1$,
at indices $q+1$ and $q$ gives
\[
 \begin{aligned}
 \Gamma_m
 ={}&3u_{m+1,q+1}-u_{m+1,q-1}
     -2u_{m+1,q+2}\\
 &+2(b_{m+1,q+2}-b_{m+1,q+1}).
 \end{aligned}
\]
For brevity, write
\[
 \rho_{m,q}
 =\frac{q(q+1)}{(m+1-q)(m+2-q)}.
\]
The binomial coefficient ratios give
\[
 \begin{aligned}
 \frac{u_{m+1,q+2}}{u_{m+1,q+1}}
 &=\frac{2(m-q)}{21(q+2)},\\
 \frac{u_{m+1,q-1}}{u_{m+1,q+1}}
 &=\frac{441}{4}\rho_{m,q},\\
 \frac{b_{m+1,q+2}-b_{m+1,q+1}}
      {u_{m+1,q+1}}
 &=\frac{N_m}{21(q+2)}(E_q+1).
 \end{aligned}
\]
Consequently,
\[
 \frac{21\Gamma_m}{u_{m+1,q+1}}
 =63-\frac{4(m-q)}{q+2}
   -\frac{9261}{4}\rho_{m,q}
   +\frac{2N_m}{q+2}(E_q+1).
\]
Using
\[
 N_m=11(m-q)-21(q+2),
\]
we can combine the terms involving $m-q$ and $N_m$ to obtain
\begin{equation}\label{eq:boundary}
 \begin{aligned}
 \Gamma_m
 &=\frac{u_{m+1,q+1}}{21}B_{m,q},\\
 B_{m,q}
 &=\frac{2N_m}{q+2}\left(E_q+\frac9{11}\right)
   +\frac{609}{11}
   -\frac{9261}{4}\rho_{m,q}.
 \end{aligned}
\end{equation}
The prefactor $u_{m+1,q+1}/21$ is positive, so
$\Gamma_m$ and $B_{m,q}$ have the same sign.

The cases $q=1,2,3$ correspond precisely to $m=5,\ldots,13$.
Direct substitution yields
\[
 \begin{array}{c|c|r|r}
 m&q&N_m&B_{m,q}\\ \hline
 5 &1&-19&-9597/20\\
 6 &1& -8&-861/4\\
 7 &1&  3&525/16\\
 8 &2&-18&-27015/16\\
 9 &2& -7&-2877/4\\
 10&2&  4&1167/5\\
 11&3&-17&-258981/40\\
 12&3& -6&-526443/220\\
 13&3&  5&147369/88
 \end{array}
\]
Thus the required signs hold in these cases.

Henceforth assume $q\ge4$.
We will use a uniform lower bound for the exponential
factor to determine the sign of $B_{m,q}$ when $N_m\ne0$.
The case $N_m=0$ will require a separate estimate.

Since
\[
 E_{q+1}=\frac{11}{2}E_q+\frac92,
\]
we have, for $q\ge0$,
\[
 \frac{E_{q+1}}{q+3}
 >
 \frac{11(q+2)}{2(q+3)}
 \frac{E_q}{q+2}
 >
 \frac{E_q}{q+2}.
\]
Therefore
\[
 \frac{2E_q}{q+2}
 \ge\frac{2E_4}{6}
 =\frac{53673}{32}
 >1000
 \qquad(q\ge4).
\]

\medskip
\noindent
\emph{Case 1: $N_m>0$.}

Since $N_m$ is an integer, $N_m\ge1$.
Moreover,
\[
 \begin{aligned}
 m+1-q
 &>\frac{21q+53}{11}>\frac{21q}{11},\\
 m+2-q
 &>\frac{21q+64}{11}>\frac{21(q+1)}{11}.
 \end{aligned}
\]
It follows that
\[
 \rho_{m,q}<\frac{121}{441},
 \qquad
 \frac{9261}{4}\rho_{m,q}
 <\frac{2541}{4}<636.
\]
The first term in $B_{m,q}$ is greater than $1000$,
and its constant term $609/11$ is positive.
Thus~\eqref{eq:boundary} gives
\[
 B_{m,q}>1000-636=364>0.
\]

\medskip
\noindent
\emph{Case 2: $N_m<0$.}

Here $N_m\le-1$, so
\[
 \frac{2N_m}{q+2}\left(E_q+\frac9{11}\right)
 <-\frac{2E_q}{q+2}
 <-1000.
\]
The last term of $B_{m,q}$ is negative, while
$609/11<56$.
Therefore
\[
 B_{m,q}<-1000+56<0.
\]

\medskip
\noindent
\emph{Case 3: $N_m=0$.}

In this case,
$b_{m+1,q+2}=b_{m+1,q+1}$, so the auxiliary contribution
to $\Gamma_m$ vanishes.
The remaining terms give
\[
 B_{m,q}
 =\frac{609}{11}-\frac{9261}{4}\rho_{m,q}.
\]
The equation $N_m=0$ gives
\[
 m+2-q=\frac{21q+64}{11}<4q
 \qquad(q\ge4).
\]
Since $m+1-q<m+2-q$, we obtain
\[
 \rho_{m,q}
 >\frac{q(q+1)}{16q^2}
 >\frac1{16}.
\]
Consequently,
\[
 B_{m,q}
 <\frac{609}{11}-\frac{9261}{64}<0.
\]

Thus $\Gamma_m>0$ when $N_m>0$, and
$\Gamma_m<0$ when $N_m\le0$.
Together with~\eqref{eq:qmstep}, this proves the assertion.
\end{proof}
	
\begin{theorem}\label{st:3.8}
For every $m\ge1$, the nonzero coefficient sequence of
$\PX{R_m}(z)$ is strictly unimodal.
For $m\ge5$, its unique mode occurs at degree
\begin{equation}\label{eq:peak}
 \mu_m=\left\lfloor\frac{149m+85}{32}\right\rfloor.
\end{equation}
\end{theorem}

\begin{proof}
For $m=1,2,3$, the nonzero coefficient sequences displayed
in the proof of Theorem~\ref{st:3.5} are strictly increasing,
and hence strictly unimodal.
The case $m=4$ follows from Example~\ref{ex:R4}.

For $m\ge5$, we prove by induction that $C_m(t)$ is strictly
unimodal with its unique mode at $q_m$.
The reversal identity will then give the corresponding
assertions for $\PX{R_m}(z)$.

For the base case, equation~\eqref{eq:cm} gives
\[
 (d_{5,0},\ldots,d_{5,6})
 =
 (21587391,27699651,13984110,3476970,358998,1712,32).
\]
Since $q_5=1$, this sequence is strictly unimodal with
its unique mode at $q_5$.

Now suppose that $C_m(t)$ is strictly unimodal with
unique mode $q=q_m$, so that
\[
 d_{m,0}<\cdots<d_{m,q}
 >d_{m,q+1}>\cdots>d_{m,m+1}.
\]
Since $q_{m+1}\in\{q,q+1\}$, we will first control the
coefficients to the left of $q$ and to the right of $q+1$,
and then compare the two coefficients at indices $q$
and $q+1$.

Retain the zero conventions from Lemma~\ref{st:3.6}.
Taking consecutive differences in~\eqref{eq:dmrecur}
gives, for $0\le k\le m+1$,
\begin{equation}\label{eq:diff}
 \begin{aligned}
 d_{m+1,k+1}-d_{m+1,k}
 ={}&21(d_{m,k+1}-d_{m,k})\\
 &+2(d_{m,k}-d_{m,k-1})\\
 &+18(b_{m,k+1}-b_{m,k}).
 \end{aligned}
\end{equation}

We also need the signs of the auxiliary differences.
For $0\le k\le m$,
\[
 b_{m,k+1}-b_{m,k}
 =\frac{b_{m,k}}{21(k+1)}
   (11m-32k-21).
\]
The prefactor is positive throughout this range.
By the definition of $q$,
\[
 -32<11m-32q-53\le0.
\]
Writing
\[
 11m-32k-21
 =(11m-32q-53)+32(q-k+1),
\]
we see that the auxiliary difference is positive for
$k<q$ and nonpositive for $k>q$ within this range.
At $k=m+1$, both auxiliary coefficients are zero.
Thus
\[
 \begin{aligned}
 b_{m,k+1}-b_{m,k}&>0
 &&(0\le k<q),\\
 b_{m,k+1}-b_{m,k}&\le0
 &&(q<k\le m+1).
 \end{aligned}
\]

\medskip
\noindent
\emph{Coefficients to the left of the candidate peak.}

For $0\le k<q$, the induction hypothesis gives
\[
 d_{m,k+1}-d_{m,k}>0,
 \qquad
 d_{m,k}-d_{m,k-1}>0.
\]
For the second inequality at $k=0$, we use
$d_{m,-1}=0$ and $d_{m,0}>0$.
The auxiliary difference is also positive.
Hence every term on the right of~\eqref{eq:diff}
is positive, and therefore
\[
 d_{m+1,0}<\cdots<d_{m+1,q}.
\]

\medskip
\noindent
\emph{Coefficients to the right of the candidate peak.}

For $q<k\le m+1$, the induction hypothesis and the zero
convention give
\[
 d_{m,k+1}-d_{m,k}<0,
 \qquad
 d_{m,k}-d_{m,k-1}<0.
\]
At the endpoint $k=m+1$, the first difference is
$-d_{m,m+1}<0$.
The auxiliary difference is nonpositive.
Thus~\eqref{eq:diff} is strictly negative, giving
\[
 d_{m+1,q+1}>d_{m+1,q+2}>\cdots>d_{m+1,m+2}.
\]

\medskip
\noindent
\emph{The remaining comparison.}

By Lemma~\ref{st:3.7},
\[
 \begin{aligned}
 d_{m+1,q}&>d_{m+1,q+1}
 &&\text{if }q_{m+1}=q,\\
 d_{m+1,q}&<d_{m+1,q+1}
 &&\text{if }q_{m+1}=q+1.
 \end{aligned}
\]
Combining the appropriate inequality with the two chains
above shows that $C_{m+1}(t)$ is strictly unimodal with
its unique mode at $q_{m+1}$.
This completes the induction.

Finally, reversal preserves strict unimodality and sends
the coefficient at degree $q_m$ in $C_m(t)$ to the
coefficient at degree $5m+1-q_m$ in $\PX{R_m}(z)$.
Thus, for $m\ge5$,
\[
 \mu_m=5m+1-q_m.
\]
Using $n-\lceil x\rceil=\lfloor n-x\rfloor$ for integers
$n$, we obtain
\[
 \begin{aligned}
 \mu_m
 &=\left\lfloor
   5m+1-\frac{11m-53}{32}
   \right\rfloor\\
 &=\left\lfloor\frac{149m+85}{32}\right\rfloor.
 \end{aligned}
\]
Together with the cases $m=1,2,3,4$, this proves
both assertions.
\end{proof}
	
\begin{remark}\label{st:3.9}
For $m\ge5$, formula~\eqref{eq:peak} gives
\[
 \frac{149m+85}{32}-1
 <\mu_m
 \le\frac{149m+85}{32}.
\]
Dividing by $|V(R_m)|=5m+1$ and letting $m\to\infty$,
we obtain
\[
 \frac{\mu_m}{|V(R_m)|}\longrightarrow\frac{149}{160}.
\]
Thus the modal degree is asymptotically $149/160$
of the number of vertices.
\end{remark}

\begin{corollary}\label{st:3.10}
For every $m\ge4$, the nonzero coefficient sequence
of $\PX{R_m}(z)$ is strictly unimodal but not log-concave.
Thus $\{R_m:m\ge4\}$ is an infinite family of trees
demonstrating that strict unimodality does not imply
log-concavity for partial Petrial polynomials.
\end{corollary}

\begin{proof}
This follows immediately from Theorems~\ref{st:3.5}
and \ref{st:3.8}.
\end{proof}

\section{Repeated path branches and eventual failure of log-concavity}

We first clarify the scope of the corank-one specialization.
Suppose that $H$ is a nonempty simple graph satisfying
\[
 \crk\A(H,F)\le1
 \qquad\text{for every }F\subseteq V(H).
\]
Then $H$ must be a path.

To see this, choose each diagonal entry to be the parity
of the degree of the corresponding vertex.
Every row then has sum zero over $\GF(2)$.
For each connected component of $H$, its indicator vector
therefore belongs to the kernel.
These vectors are linearly independent, so the corank
is at least the number of components.
The assumed upper bound of one forces $H$ to be connected.

Every simple graph admits a nonsingular diagonal
perturbation.
We prove this by induction on the number of vertices,
starting with the empty matrix, whose determinant is $1$.
For a nonempty graph $G$, choose a vertex $v$ and,
by induction, fix a nonsingular diagonal perturbation
on $G-v$.
With these diagonal entries fixed, the determinant of
the full matrix has the form
\[
 \delta_v+\beta,
 \qquad \beta\in\GF(2),
\]
where $\beta$ is independent of $\delta_v$.
Indeed, the coefficient of $\delta_v$ is the determinant
of the chosen matrix on $G-v$, which equals $1$.
Exactly one choice of $\delta_v$ therefore makes the
full matrix nonsingular.

Thus $H$ has both singular and nonsingular diagonal
perturbations.
Together with the corank bound, this implies that
$C_H(t)$ has degree $1$ and a positive constant term.
Consequently, both $C_H(t)$ and $\PX{H}(z)$ have exactly
two nonzero terms.
The binomial characterization for connected circle graphs
was established by Feng, Yan and Zheng~\cite{FengYanZheng}
and extended to all connected simple graphs by Deng and
Zhang~\cite{DengZhang}.
Applying the latter characterization, we conclude that
$H$ is a path.

The general attachment identity~\eqref{eq:root-general}
applies to arbitrary rooted simple graphs.
The specialization~\eqref{eq:rootjoin} is therefore
restricted to rooted paths.
We next give a criterion for eventual failure of
log-concavity in this setting.

\begin{theorem}\label{st:3.11}
Let $(H,c)$ satisfy the hypothesis of
Corollary~\ref{st:3.1}, with perturbation counts
$a,b,\gamma$.
If $b>0$ and $\gamma>0$, then there exists a positive
integer $M$ such that, for every $m\ge M$, the nonzero
coefficient sequence of $\PX{T_m(H,c)}(z)$ is not
log-concave.
\end{theorem}

\begin{proof}
Fix $(H,c)$, so that $a,b,\gamma$ remain constant as
$m$ varies.
Put
\[
 D_m(t)=C_{T_m(H,c)}(t),\qquad
 \alpha=\frac ab,\qquad
 \rho=1+\frac{\gamma}{b}>1.
\]
For $m\ge1$, let $e_{m,0},e_{m,1},e_{m,2}$ denote
the three highest-degree coefficients of $D_m(t)$,
corresponding to $t^{m+1},t^m,t^{m-1}$, respectively.

Extracting these coefficients from~\eqref{eq:rootjoin}
gives
\[
 \begin{aligned}
 e_{m,0}
 &=b^m,\\
 e_{m,1}
 &=b^m(1+m\alpha),\\
 e_{m,2}
 &=b^m\left[
   m\alpha+\binom m2\alpha^2+2(\rho^m-1)
   \right].
 \end{aligned}
\]
These are three consecutive positive coefficients.
To prove failure of log-concavity, it suffices to show that
\[
 \Delta_m:=e_{m,1}^2-e_{m,0}e_{m,2}<0.
\]
After dividing by the positive factor $b^{2m}$, we obtain
\begin{equation}\label{eq:general-gap}
 \begin{aligned}
 \frac{\Delta_m}{b^{2m}}
 &=(1+m\alpha)^2
   -m\alpha-\binom m2\alpha^2-2(\rho^m-1)\\
 &=3+m\alpha+\frac{m(m+1)}2\alpha^2-2\rho^m.
 \end{aligned}
\end{equation}
Since $\rho>1$,
\[
 \frac{3+m\alpha+\frac{m(m+1)}2\alpha^2}{\rho^m}
 \longrightarrow0.
\]
In particular, this ratio is less than $2$ for all
sufficiently large $m$.
Equation~\eqref{eq:general-gap} then gives $\Delta_m<0$.

Thus the nonzero coefficient sequence of $D_m(t)$
is not log-concave for all sufficiently large $m$.
Reversal preserves log-concavity, so the same conclusion
holds for $\PX{T_m(H,c)}(z)$.
\end{proof}

\begin{remark}\label{rem:tail-defect}
Under the assumptions of Theorem~\ref{st:3.11}, retain
the notation
\[
 \alpha=\frac ab,\qquad
 \rho=1+\frac{\gamma}{b}>1.
\]
Let $e_{m,0},e_{m,1},e_{m,2}$ be the three highest-degree
coefficients of $D_m(t)=C_{T_m(H,c)}(t)$, listed in
decreasing degree order.
Then
\[
 \frac{e_{m,1}^{\,2}}{e_{m,0}e_{m,2}}
 =
 \frac{(1+m\alpha)^2}
 {m\alpha+\binom m2\alpha^2+2(\rho^m-1)}
 \longrightarrow0.
\]
In the fraction on the right, the numerator is $O(m^2)$,
while the denominator is asymptotic to $2\rho^m$,
which proves the limit.
Since log-concavity requires this ratio to be at least $1$,
the limit quantifies the failure at this triple.
Other log-concavity inequalities may also fail.
\end{remark}

We now apply the criterion to paths rooted at internal
vertices.
Let
\[
 J_0=0,\qquad J_1=1,\qquad
 J_{r+1}=J_r+2J_{r-1}\quad(r\ge1)
\]
be the Jacobsthal numbers.
Their closed form is
\[
 J_r=\frac{2^r-(-1)^r}{3}
 \qquad(r\ge0).
\]
We use the convention that $P_0$ is empty and
$C_{P_0}(t)=1$.
By~\eqref{eq:leaf} and coefficient reversal,
\[
 C_{P_s}(t)
 =C_{P_{s-1}}(t)+2C_{P_{s-2}}(t)
 \qquad(s\ge2).
\]
Together with $C_{P_1}(t)=1+t$, this recurrence gives
\[
 C_{P_s}(t)=J_{s+1}+J_st
 \qquad(s\ge0).
\]
We also record the addition identity
\[
 J_{i+j+1}
 =J_{i+1}J_{j+1}+2J_iJ_j
 \qquad(i,j\ge0),
\]
which follows directly from the closed form.

\begin{corollary}\label{st:3.12}
Let $N\ge3$ and let $c$ be an internal vertex of $P_N$.
Then there exists a positive integer $M$ such that,
for every $m\ge M$, the nonzero coefficient sequence
of $\PX{T_m(P_N,c)}(z)$ is not log-concave.
\end{corollary}

\begin{proof}
Write
\[
 P_N-c=P_i\sqcup P_j,
 \qquad i,j\ge1,\qquad i+j=N-1.
\]
By the path formula, every diagonal perturbation of $P_N$
has corank at most $1$, and
\[
 a=J_{N+1},\qquad b+\gamma=J_N.
\]
We determine $\gamma$ by examining the two arms.

A nonzero kernel vector of a perturbed path matrix has
coordinate $1$ at each endpoint.
Indeed, if its coordinate at one endpoint were zero,
the successive row equations would force all coordinates
to vanish.

If either arm matrix is singular, the kernel argument
in the proof of Lemma~\ref{st:root-general} therefore
forces the coordinate at $c$ to be zero.
Such a perturbation cannot be counted by $\gamma$.

Conversely, suppose both arm matrices are nonsingular.
Setting the coordinate at $c$ to $1$ uniquely determines
the coordinates on both arms.
The central equation then determines exactly one choice
of the diagonal entry at $c$.
Since the full path matrix has corank at most one,
the resulting perturbation is counted by $\gamma$.

There are $J_{i+1}$ nonsingular perturbations of $P_i$
and $J_{j+1}$ of $P_j$.
Hence
\[
 \gamma=J_{i+1}J_{j+1}.
\]
By the Jacobsthal addition identity,
\[b=J_N-\gamma=J_{i+j+1}-J_{i+1}J_{j+1}=2J_iJ_j.
\]
All Jacobsthal numbers with positive index are positive.
Since $i,j\ge1$, we have $b>0$ and $\gamma>0$.
The conclusion now follows from Theorem~\ref{st:3.11}.
\end{proof}

\begin{remark}\label{st:3.13}
Let $\ell_G$ and $h_G$ denote the coefficients of the
lowest-degree and highest-degree nonzero terms of
$\PX{G}(z)$, respectively.
By coefficient reversal and~\eqref{eq:explicit},
\[
 \begin{aligned}
 \ell_{R_m}
 &=d_{m,m+1}=2^m,\\
 h_{R_m}
 &=d_{m,0}
 =21^{m-1}(21+18m).
 \end{aligned}
\]
Thus $\ell_{R_m}<h_{R_m}$ for every $m\ge1$, in agreement
with the extremal coefficient inequality
of Feng, Guo and Yan~\cite{FengGuoYan}.
However, for every $m\ge4$, the nonzero coefficient
sequence of $\PX{R_m}(z)$ is not log-concave by
Theorem~\ref{st:3.5}.
Hence the inequality $\ell_G<h_G$ alone does not
guarantee log-concavity, even when $G$ is a tree.
\end{remark}

\section{Concluding remarks}
\label{sec:conclusions}

We constructed a family of trees $R_m$ whose partial
Petrial polynomials are strictly unimodal for every
$m\ge1$, but log-concave exactly when $1\le m\le3$.
We also determined their exact modes for $m\ge5$.

More generally, our rooted attachment formula shows
that repeated attachment of any path rooted at an
internal vertex leads to failure of log-concavity
for all sufficiently large multiplicities.
These results motivate the following question
about unimodality.

\begin{problem}[Unimodality]
Is the nonzero coefficient sequence of $\PX{T}(z)$
unimodal for every tree $T$?
More generally, does unimodality hold for every
simple graph?
\end{problem}

Known formulas also provide immediate examples of
log-concavity.
For instance, the path formula of Yan and Li~\cite{YanLi}
is
\[
 \PX{P_n}(z)
 =J_nz^{n-1}+J_{n+1}z^n
 \qquad(n\ge1).
\]
Its two positive coefficients form a log-concave
sequence, since there is no internal inequality to check.
Together with the families constructed here, these
examples motivate the following structural question.

\begin{problem}[Structural log-concavity]
Which trees have log-concave partial Petrial polynomials?
In particular, does every tree of maximum degree at
most $3$ have this property?
\end{problem}

\section*{Acknowledgements} 
This work is supported by NSFC (No. 12471326).

\end{document}